\documentclass[11pt,reqno]{amsart}
\usepackage{amsmath,amsthm,amssymb,verbatim,enumerate,ifthen}
\usepackage[mathscr]{eucal}
\usepackage[utf8]{inputenc}
\usepackage[T1]{fontenc}
\def\N{\mathbb{N}}
\def\R{\mathbb{R}}

\def\F{\mathscr{F}}
\def\K{\mathscr{K}}
\def\TB{\mathscr{T}}

\def\B{\mathscr{B}}

\def\diam{\mathop{\mbox{\rm diam}}}
\let\eps\varepsilon

\newtheorem{theorem}{Theorem}
\newtheorem*{theorem*}{Theorem}
\def\Thm#1#2{\ifthenelse{\equal{#1}{*}}{\begin{theorem*}#2\end{theorem*}}
  {\begin{theorem}\label{T#1}#2\end{theorem}}}
\newtheorem{Atheorem}{Theorem}

\def\thm#1{Theorem~\ref{T#1}}

\newtheorem{proposition}[theorem]{Proposition}
\newtheorem*{proposition*}{Proposition}
\def\Prp#1#2{\ifthenelse{\equal{#1}{*}}{\begin{proposition*}#2\end{proposition*}}
             {\begin{proposition}\label{P#1}#2\end{proposition}}}

\newtheorem{corollary}[theorem]{Corollary}
\newtheorem*{corollary*}{Corollary}
\def\Cor#1#2{\ifthenelse{\equal{#1}{*}}{\begin{corollary*}#2\end{corollary*}}
             {\begin{corollary}\label{C#1}#2\end{corollary}}}

\newtheorem{lemma}[theorem]{Lemma}
\newtheorem*{lemma*}{Lemma}
\def\Lem#1#2{\ifthenelse{\equal{#1}{*}}{\begin{lemma*}#2\end{lemma*}}
             {\begin{lemma}\label{L#1}#2\end{lemma}}}
\def\lem#1{Lemma~\ref{L#1}}

\newtheorem{example}[theorem]{Example}
\newtheorem*{example*}{Example}
\def\Exa#1#2{\ifthenelse{\equal{#1}{*}}{\begin{example*}\rm #2\end{example*}}
             {\begin{example}\label{Ex#1}\rm #2\end{example}}}

\newtheorem{problem}[theorem]{Problem}

\theoremstyle{definition}
\newtheorem{definition}[theorem]{Definition}

\newtheorem{remark}[theorem]{Remark}
\newtheorem*{remark*}{Remark}
\def\Rem#1#2{\ifthenelse{\equal{#1}{*}}{\begin{remark*}\rm #2\end{remark*}}
             {\begin{remark}\label{R#1}\rm #2\end{remark}}}

\def\eq#1{{\rm(\ref{E#1})}}
\def\Eq#1#2{\ifthenelse{\equal{#1}{*}}
  {\begin{equation*}\begin{aligned}#2\end{aligned}\end{equation*}}
  {\begin{equation}\begin{aligned}\label{E#1}#2\end{aligned}\end{equation}}}

\long\def\comment#1{}

\begin{document}
\begin{flushright}
\end{flushright}
\vspace{5mm}

\date{\today}
\title{Hutchinson's Theorem in Semimetric Spaces}

\author[Zs. P\'ales]{Zsolt P\'ales}
\address[Zs. P\'ales]{Institute of Mathematics, University of Debrecen, H-4002 Debrecen, Pf.\ 400, Hungary}
\email{pales@science.unideb.hu}

\author[M. Kocsis]{Mátyás Kocsis}
\address[M. Kocsis]{Faculty of Science and Technology, University of Debrecen, H-4002 Debrecen, Pf.\ 400, Hungary}
\email{kocsis.matyas21@gmail.com}

\subjclass[2010]{54E25,28A80}
\keywords{semimetric space, Hausdorff--Pompeiu distance, comparison function, contraction, fixed point theorem, fractal, stability}


\thanks{The research of the second author was supported by the K-134191 NKFIH Grant and by the EFOP-3.6.1-16-2016-00022 project, which was co-financed by the European Union and the European Social Fund.}

\begin{abstract}
One of the important consequences of the Banach Fixed Point Theorem is Hutchinson's theorem which states the existence and uniqueness of fractals in complete metric spaces. The aim of this paper is to extend this theorem for semimetric spaces using the results of Bessenyei and Páles published in 2017. In doing so, some properties of semimetric spaces as well as of the fractal space are investigated. We extend Hausdorff's theorem to characterize compactness and Blaschke's theorems to characterize the completeness of the fractal space. Based on these preliminaries, an analogue of Hutchinson's Theorem in the setting of semimetric spaces is proved and finally, error estimates and stability of fractals are established as well.
\end{abstract}

\maketitle

\section{Introduction}

Let $X$ be a nonempty set, and let $T_1,\dots,T_n\colon X\to X$. A nonempty subset $F\subset X$ is called a \emph{fractal} with respect to the system $(T_1,\dots,T_n)$ if it fulfills the so-called \emph{invariance equation}:
\Eq{A}{
 F=\bigcup_{k=1}^nT_k(F).
}
It is clear that if the right hand side of this equality is considered as set-to-set mapping, then $F$ is its fixed point. Therefore, the main question is to clarify what conditions for $X$ and for the system $(T_1,\ldots,T_n)$ ensure that $F$ exists and be unique. The theory of fixed points provides a rich toolkit to study such questions. We only mention some standard and important monographs: Berinde \cite{Ber07}, Granas and Dugundji \cite{GraDug03}, Rus \cite{Rus01}, Rus, Petruşel and Petruşel \cite{RusPetPet08}, and Zeidler \cite{Zei86}. In these books the theory of fixed points is mainly treated in complete metric and Banach spaces. In his seminal paper \cite{Hut81a} Hutchinson used the Banach Contraction Principle (cf.\ Banach \cite{Ban22} and Cacciopoli \cite{Cac30}) to establish the existence of fractals assuming that $X$ is a complete metric space and $T_1,\dots,T_n$ are contractions. There are only a few results in the much more general setting of semimetric spaces. In the paper \cite{BesPal17}, fixed points theorems in complete semimetric spaces were established by Bessenyei and Páles which generalized the basic results of Browder \cite{Bro68,Bro79} and Matkowski \cite{Mat75} and which enjoy many extensions (cf.\ Miculescu and Mihail \cite{MicMih17}, Mitrović and Hussain \cite{MitHus19}). In this paper we aim to apply the results of the paper \cite{BesPal17} to prove the existence and uniqueness of fractals.

Throughout this paper, let $(X,d)$ be a semimetric space, which means that $X$ is a nonempty set and $d:X\times X\to \R$ is a nonnegative-valued symmetric function which vanishes exactly at the diagonal points of the Cartesian product $X\times X$. The theory of semimetric spaces was investigated in the the last century in the papers Burke \cite{Bur72}, Galvin and Shore \cite{GalSho84}, McAuley \cite{McA56} and Wilson \cite{Wil31b}.  Some important recent developments are due to Chrąszcz, Jachymski and Turoboś \cite{ChrJacTur18,ChrJacTur19,JacTur20} and Dung and Hang \cite{VanHan17}. Recalling some of these results, in Section 2, we describe the most important definitions (such as convergence, completeness, boundedness, topology, etc.) and basic results about semimetric spaces. The main result of this section is the extension of Hausdorff's Theorem about the characterization of compactness in metric spaces to the semimetric setting.

In terms of a semimetric $d$, we define the Hausdorff--Pompeiu distance of two nonempty subsets analogously to the standard definition (cf.\ \cite{Hau49}, \cite{Pom05}) in metric spaces in Section 3. We establish some of its basic properties and we show that the class of nonempty bounded an closed sets and also the class of compact nonempty sets forms a semimetric space equipped with the Hausdorff--Pompeiu distance. We characterize the completeness of these spaces and thus we extend Blaschke's celebrated results (cf.\ \cite{Bla49}) to the semimetric setting. 

The main goals of this paper are reached in Section 4, where we recall the notion of a comparison function (\cite{Bes15}, \cite{Bes16}, \cite{BesPal17}, \cite{JacMatSwi95}) and we define the concept of contractions with respect to comparison functions. First we state the results of Bessenyei and Páles \cite{BesPal17} which establish the existence, stability of fixed points of contractions as well as we prove error estimates for the so-called Picard iteration. Using these results, we obtain a generalization of Hutchinson's above described result. We also provide  error estimates for the corresponding iteration and prove the stability of fractals with respect to pointwise convergence of contractions. We mentions here that in order to obtain such results, one may not need a generalization of Blaschke's theorem. An approach where the use of Blaschke's theorem is avoided, has been presented by Bessenyei and Pénzes \cite{BesPen21}.

\section{Basic terminology and results about semimetric spaces}

Throughout this paper, let $\R_+$ and $\overline{\R}_+$ denote the sets of nonnegative real numbers and nonnegative extended real numbers, respectively.
We say that $\Phi\colon\R_+^2\to\overline{\R}_+$ is a \emph{triangle function} for the semimetric $d$ (cf.\ \cite{BesPal17}), if $\Phi$ is symmetric and monotone increasing in both of its arguments, $\Phi(0,0)=0$ and, for all $x,y,z\in X$, the \emph{triangle inequality}
\Eq{*}{
 d(x,y)\le\Phi\bigl(d(x,z),d(z,y)\bigr)
}
holds. For a semimetric space $(X,d)$, define the function $\Phi_d\colon\R^2_+\to\overline{\R}_+$ by
\Eq{basic}{
 \Phi_d(u,v):=\sup\{d(x,y)\mid\exists p\in X: d(p,x)\le u,\,d(p,y)\le v\}.
}
A simple and direct calculation shows that $\Phi_d$ is a triangle function for $d$. This function is called the \emph{basic triangle function related to $d$}. It is easy to see that the basic triangle function is optimal, that is, if $\Phi$ is a triangle function for $d$, then $\Phi_d\le\Phi$ holds.

A triangle function $\Phi$ is called \emph{regular} if it is continuous at $(0,0)$. A semimetric space $(X,d)$ is called \emph{regular} if it admits a regular triangle function. Clearly, in a regular semimetric space, the basic triangle function is regular.

The notions of a \emph{convergent sequence} and a \emph{Cauchy sequence} in a semimetric space are defined in the standard way \cite{Bur72}. A semimetric space is termed to be \emph{complete}, if each Cauchy sequence of the space is convergent \cite{GalSho84}. Under the \emph{diameter} of a set in a semimetric space, we mean the supremum of distances taken over the pairs of points of the set. A set is called \emph{bounded} if its diameter is finite.

Let $(X,d)$ be a semimetric space. Then the \emph{ball of radius $r\in(0,\infty]$ centered at $p\in X$} is defined as 
\Eq{*}{
\B(p,r):=\{x\in X\mid d(x,p)<r\}.
}
The following lemma from \cite{BesPal17} characterizes regular semimetric spaces and summarizes their basic properties.

\Lem{regsem}{A semimetric space $(X,d)$ is regular if and only if
\Eq{diameter}{
 \lim_{r\to 0}\sup_{p\in X}\diam \B(p,r)=0.
}
Furthermore, in a regular semimetric space, convergent sequences have a unique limit and possess the Cauchy property.}

In view of this lemma, convergence and Cauchy property are equivalent to each other in complete and regular semimetric spaces. 

We define the \emph{topology} of a semimetric space $(X,d)$ in terms of interior points. We say that a subset $A\subseteq X$ is \emph{open} if every element $p\in A$ is also an \emph{interior point}, i.e., there is an $r>0$ such that $\B(p,r)\subseteq A$. It is easy to see that the open sets form a topology. However, the balls in a semimetric space  are not necessarily open.

\Lem{T}{Let $(X,d)$ be a semimetric space. Let $A\subseteq X$ be a closed set and $(x_n)$ be a convergent sequence of elements from $A$. Then the limit of $(x_n)$ also belongs to $A$. If the semimetric space is also regular, then, for every $\eps>0$, there exists $r>0$ such that, for all $x\in X$, the inclusion $\B(x,r)\subseteq \B(x,\eps)^\circ$ holds, consequently, the topology is Hausdorff.}

\begin{proof}
Let $p$ be the limit of a sequence $(x_n)$ belonging to $A$. If $p\not\in A$, then $p\in X\setminus A$, which is open, and hence, for some $r>0$, we have that $\B(p,r)\subseteq X\setminus A$. Therefore, $\B(p,r)\cap A=\emptyset$. On the other hand, by the convergence $x_n\to p$, there exists $n_0\in\N$ such that for all $n\geq n_0$, we have $d(x_n,p)<r$. This implies $x_n\in\B(p,r)\cap A$, which is a contraction.

Let $\eps>0$. By the regularity of the semimetric space, there exists $r>0$ such that $\Phi_d(r,r)<\eps$. Let $x\in X$, $y\in\B(x,r)$ and $z\in\B(y,r)$. Then
\Eq{*}{
  d(z,x)\leq\Phi_d(d(z,y),d(y,x))\leq\Phi_d(r,r)<\eps.
}
Hence, $z\in\B(x,\eps)$, which implies $\B(y,r)\subseteq\B(x,\eps)$. This proves that $y$ is an interior point of $\B(x,\eps)$, therefore $\B(x,r)\subseteq \B(x,\eps)^\circ$.

Let $x,y$ be distinct points of $X$. Then there exists $\eps>0$ such that $\Phi_d(\eps,\eps)<d(x,y)$. This implies that the balls $\B(x,\eps)$ and $\B(y,\eps)$ are disjoint. By the previous statement, the sets $\B(x,\eps)^\circ$ and $\B(y,\eps)^\circ$ are disjoint open sets containing $x$ and $y$, respectively, which proves the Hausdorff property.
\end{proof}

If $(X,d)$ and $(Y,\rho)$ are semimetric spaces, then a mapping $T:X\to Y$ is called \emph{Lipschitzian} if there exists $L\geq0$ such that, for all $x,y\in X$,
\Eq{Lip}{
  \rho(T(x),T(y))\leq L d(x,y).
}

\Lem{cont}{Let $(X,d)$ and $(Y,\rho)$ be semimetric spaces and let $T:X\to Y$ be a Lipschitzian mapping. Then $T$ is continuous with respect to the topology of the semimetric spaces (i.e., the inverse image of any open subset of $Y$ by $T$ is open in $X$).}

\begin{proof} Assume that \eq{Lip} holds, for all $x,y\in X$, with some constant $L\geq0$. To show the continuity of $T$, let $U\subseteq Y$ be an open set and let $p\in T^{-1}(U)$. Then, $T(p)\in U$. Since $U$ is open, we have that $T(p)$ is an interior point to $U$, i.e., there exists $r>0$ such that $\B(T(p),r)\subseteq U$. We show that $\B(p,r/L)\subseteq T^{-1}(U)$. Indeed, if
$x\in\B(p,r/L)$, then $d(x,p)<r/L$. Hence, in view of \eq{Lip}, we have that $\rho(T(x),T(p))\leq Ld(x,p)<r$, which proves that $T(x)\in \B(T(p),r)\subseteq U$. Hence $x\in T^{-1}(U)$, showing that $\B(p,r/L)\subseteq T^{-1}(U)$. Thus we have proved that every point of $T^{-1}(U)$ is an interior point to this set. Therefore, $T^{-1}(U)$ has to be open and, consequently, $T$ is continuous.
\end{proof}

A subset $H$ of a semimetric space $(X,d)$ is called \emph{totally bounded} if, for all $\eps>0$, there exists a finite subset $P:=\{p_1,\dots,p_n\}\subseteq X$, such that $H\subseteq\bigcup_{i=1}^n\B(p_i,\eps)$. If this holds, then $P$ will be called an $\eps$-net for $H$. It is clear that totally bounded sets are always bounded. The family of all nonempty closed and totally bounded subsets of the semimetric space $(X,d)$ will be denoted by $\TB(X)$.

A subset $H$ of a semimetric space $(X,d)$ is called \emph{compact} if every open cover of $H$ contains a finite subcover of $H$. The family of all nonempty compact subsets of the semimetric space $(X,d)$ will be denoted by $\K(X)$. 

The following result extends Hausdorff's Theorem about the characterizations of compact sets in metric spaces to the semimetric space setting.

\Thm{Haus}{Let $(X,d)$ be a complete regular semimetric space. Then a subset of $X$ is compact if and only if it is closed and totally bounded.}

\begin{proof}
Assume first that $H\subseteq X$ is a compact subset of $X$. 

To prove the total boundedness of $H$, let $\eps>0$ be arbitrary. Then
$\{\B(p,\eps)^\circ:p\in H\}$ is an open cover of $H$. By the compactness of $H$, we have that 
\Eq{*}{
H\subseteq\bigcup_{i=1}^n \B(p_i,\eps)^\circ
\subseteq\bigcup_{i=1}^n \B(p_i,\eps)
}
for some finite set $\{p_1,\dots,p_n\}\subseteq H$. This shows that $H$ is totally bounded.

To verify that $H$ is closed, let $p\in X\setminus H$. We are going to show that $p$ is an interior point of $X\setminus H$. Define, for $r>0$,
\Eq{*}{
  U_r:=\{x\in X\colon d(x,p)>\Phi(r,r)\},\qquad
  V_r:=\{x\in X\colon d(x,p)>r\}.
}
Then, for $x\in U_r$ and $y\in\B(x,r)$, we have
\Eq{*}{
  \Phi(r,r)<d(x,p)\leq\Phi(d(x,y),d(y,p))\leq\Phi(r,d(y,p)),
}
which implies that $r<d(y,p)$. Therefore, $y\in V_r$, which yields that $\B(x,r)\subseteq V_r$. This shows that $x$ is an interior point of $V_r$, hence $U_r\subseteq V_r^\circ$. 

Observe that the family $\{U_r\colon r>0\}$ covers $H$. Indeed, if $x\in H$, then $d(x,p)>0=\Phi(0,0)$. By the regularity of the triangle function, we can find a positive number $r$ such that $d(x,p)>\Phi(r,r)$, which yields that $x\in U_r$. Then, by the inclusion $U_r\subseteq V_r^\circ$, the family the family $\{V_r^\circ\colon r>0\}$ is an open cover for $H$. By the compactness, we can find $r_1,\dots,r_n>0$ such that $\{V_{r_i}^\circ\colon i\in\{1,\dots,n\}\}$ also covers $H$. This implies that $\{V_{r_i}\colon i\in\{1,\dots,n\}\}$ is a covering of $H$. 
It is obvious that $\bigcup_{i=1}^n V_{r_i}=V_{r_0}$, where $r_0:=\min(r_1,\dots,r_n)$, hence $V_{r_0}$ contains $H$.
Taking the complements of each set, we get that
\Eq{*}{
 \B(p,r)\subseteq X\setminus V_{r_0}\subseteq X\setminus H.
}
Consequently, $p$ is an interior point to $X\setminus H$, which demonstrates that $H$ is closed.

To prove the the sufficiency, let $H$ be a closed and totally bounded subset of $X$. Assume that $H$ is not compact, that is, there exists an open cover $\{C_i\}_{i\in I}$ of $H$ such that it has no finite subcover of $H$. First, construct a decreasing sequence $(\eps_n)$ of positive numbers converging to $0$ and  satisfying the following two properties:
\Eq{*}{
  \eps_1:=1,\qquad 
  \Phi(\eps_{k+1},\eps_{k+1})<\eps_k \qquad(k\in\N).
}
Choose a finite $\eps_1$-net, denoted as $P_1$, for $H_0:=H$. Then there exists an element $p_1\in P_1$ such that $\{C_i\}_{i\in I}$ has no finite subcover of $H_1:=\B(p_1,\eps_1)\cap H_0$. After this, let $P_2$ be a finite $\eps_2$-net for $H_1$. Then there exists an element $p_2\in P_2$ such that $H_2:=\B(p_2,\eps_2)\cap H_1$ has no finite subcover from $\{C_i\}_{i\in I}$. We proceed recursively. Assume that we have constructed $H_{n-1}:=\B(p_{n-1},\eps_{n-1})\cap H_{n-2}$ such that $H_{n-1}$ cannot be covered by any finite subfamily of $\{C_i\}_{i\in I}$. Let $P_n$ be a finite $\eps_n$-net for $H_{n-1}$. Then, for some $p_n\in P_n$, the intersection $H_n:=\B(p_n,\eps_n)\cap H_{n-1}$ cannot be covered by a finite subcover of $\{C_i\}_{i\in I}$.

We are now going to show that $(p_n)$ is a Cauchy sequence. The set $H_n$ cannot be empty, therefore there exists an element $q_n\in H_n$. Then $q_n\in\B(p_{n-1},\eps_{n-1})\cap\B(p_n,\eps_n)$, therefore, for all $n>2$,
\Eq{en}{
  d(p_{n-1},p_n)\leq\Phi(d(p_{n-1},q_n), d(q_n,p_n))\leq\Phi(\eps_{n-1},\eps_n)\leq\Phi(\eps_{n-1},\eps_{n-1})\leq\eps_{n-2}.
}
Using induction on $k$, we show that, for all $k\in\N$ and $m>2$,
\Eq{ek}{
  d(p_m,p_{m+k})\leq \eps_{m-2}.
}
If $k=1$, then the statement follows from \eq{en} with $n=m+1$. Assume that \eq{ek} has been proved for some $k$. Then, by applying \eq{en} and the inductive hypothesis, we get
\Eq{*}{
  d(p_m,p_{m+(k+1)})
  \leq\Phi(d(p_m,p_{m+1}),d(p_{m+1},p_{(m+1)+k}))
  \leq\Phi(\eps_{m-1},\eps_{m-1})\leq \eps_{m-2}.
}
This shows that \eq{ek} also hold for $k+1$ instead of $k$.
The sequence $(\eps_n)$ being a nullsequence, the inequality \eq{ek} implies that $(p_n)$ is a Cauchy sequence. 

By the completeness of the semimetric space, $p_n$ converges to some element $p\in X$. We have that
\Eq{*}{
  d(q_n,p)\leq\Phi(d(q_n,p_n),d(p_n,p))
  \leq \Phi(\eps_n,d(p_n,p)),
}
which shows that $(d(q_n,p))$ is a nullsequence, and hence, $(q_n)$ converges to $p$. On the other hand, for all $n\in\N$, we have that $q_n\in H_n\subseteq H$, therefore, the closedness of $H$ implies that $p$ must belong to $H$. Thus, on of the elements of the open cover $\{C_i\}_{i\in I}$ of $H$, should contain $p$, say we have $p\in C_\alpha$ for some $\alpha\in I$. Due to the openness of $C_\alpha$, there exists $r>0$ such that $\B(p,r)\subseteq C_\alpha$. Choose $\rho>0$ so that $\Phi(\rho,\rho)<r$. Then there exists $n\in\N$ such that $\eps_n<\rho$ and $d(p_n,p)<\rho$. For all $x\in H_n$, we have that $x\in\B(p_n,\eps_n)$, whence we obtain
\Eq{*}{
  d(x,p)
  \leq \Phi(d(x,p_n),d(p_n,p))
  \leq \Phi(\eps_n,d(p_n,p))
  \leq \Phi(\rho,\rho)<r.
}
This proves that $H_n\subseteq \B(p,r)\subseteq C_\alpha$, which contradicts the property that $H_n$ cannot be covered by any finite subsystem of $\{C_i\}_{i\in I}$. The contradiction so obtained shows that $H$ can be covered by a finite subsystem of $\{C_i\}_{i\in I}$.
\end{proof}

\section{The semimetric spaces of nonempty closed and bounded and nonempty compact sets}

Let $(X,d)$ be semimetric space. For two nonempty subsets $A,B\subseteq X$, their \emph{Hausdorff--Pompeiu distance} is defined by
\Eq{*}{
 D(A,B):=\inf\Bigl\{\eps\in(0,\infty]\mid 
  A\subseteq\bigcup_{b\in B}\B(b,\eps),\,B\subseteq\bigcup_{a\in A}\B(a,\eps)\Bigr\}.
}

We note that the Hausdorff--Pompeiu distance of two subsets can be infinite.
The following useful statement is an immediate consequence of the definition.

\Lem{AB}{Let $(X,d)$ be a semimetric space and let $A_1,\dots,A_n,B_1,\dots,B_n\subseteq X$. Then
\Eq{AB}{
  D(A_1\cup\dots\cup A_n,B_1\cup\dots\cup B_n)
  \leq\max\big(D(A_1,B_1),\dots,D(A_n,B_n)\big).
}}

\begin{proof} If any of the distances $D(A_i,B_i)$ is infinite then the statement is trivial. Therefore, we may assume that $\max\big(D(A_1,B_1),\dots,D(A_n,B_n)\big)<\infty$. Let $\eps>\max\big(D(A_1,B_1),\dots,D(A_n,B_n)\big)$ be arbitrary. Then, for all $j\in\{1,\dots,n\}$, 
\Eq{*}{
  A_{j}\subseteq \bigcup_{x\in B_{j}}\B(x,\eps) \qquad\mbox{and}\qquad
  B_{j}\subseteq \bigcup_{x\in A_{j}}\B(x,\eps). 
}
It follows from these inclusions that, for all $j\in\{1,\dots,n\}$,
\Eq{*}{
  A_{j}\subseteq \bigcup_{x\in B_1\cup\dots\cup B_n}\B(x,\eps) \qquad\mbox{and}\qquad
  B_{j}\subseteq \bigcup_{x\in A_1\cup\dots\cup A_n}\B(x,\eps),
}
whence we get that
\Eq{*}{
  A_1\cup\dots\cup A_n\subseteq \bigcup_{x\in B_1\cup\dots\cup B_n}\B(x,\eps) \qquad\mbox{and}\qquad
  B_1\cup\dots\cup B_n\subseteq \bigcup_{x\in A_1\cup\dots\cup A_n}\B(x,\eps).
}
These inclusions yield that $D(A_1\cup\dots\cup A_n,B_1\cup\dots\cup B_n)\leq\eps$. Upon taking the limit $$\eps\longrightarrow\max\big(D(A_1,B_1),\dots,D(A_n,B_n)\big),$$ we obtain that \eq{AB} holds.
\end{proof}

We say that a sequence of subsets \emph{$H_k\subseteq X$ converges to $H_0$ with respect to $D$} if the sequence $(D(H_k,H_0))$ tends to $0$ as $k\to +\infty$. 

\Lem{H_k}{Let $(X,d)$ be a semimetric space and let $(H_{k,1}),\dots,(H_{k,n})$ be sequences of sets that converge to some sets $H_{0,1},\dots,H_{0,n}$ with respect to $D$, respectively. For $k\in\N\cup\{0\}$, define the set $H_k$ by
\Eq{*}{
 H_k:=H_{k,1}\cup\dots\cup H_{k,n}.
}
Then the sequence $(H_k)$ converges to $H_0$ with respect to $D$.}

\begin{proof} Applying the previous lemma, we have that
\Eq{*}{
  D(H_k,H_0)=D\big(H_{k,1}\cup\dots\cup H_{k,n},H_{0,1}\cup\dots\cup H_{0,n}\big)
  \leq\max\big(D(H_{k,1},H_{0,1}),\dots,D(H_{k,n},H_{0,n})\big).
}
By the assumptions of this lemma, the right-hand side of this inequality tends to zero and hence the sequence $(D(H_k,H_0))$ also converges to zero.
\end{proof}

In the sequel, let $\F(X)$ denote the family of all nonempty, closed and bounded subsets of $X$. In what follows, we show that this class of sets forms a semimetric space with $D$ under certain conditions.

\Thm{metric}{Let $(X,d)$ be a semimetric space with an upper semicontinuous triangle function $\Phi\colon\R_+^2\to\R_+$ for $d$. Then $(\F(X),D)$ is also a semimetric space and $\Phi$ is also a triangle function for $D$.}

\begin{proof}
First we point out that $D$ has finite values. Let $A,B\in\F(X)$ be fixed. Then, by their boundedness, $\alpha:=\diam(A)$ and $\beta:=\diam(B)$ are finite. Let $a_0\in A$ and $b_0\in B$ be fixed. Then, for every $a\in A$ and $b\in B$, we have 
\Eq{*}{
d(a,b)&\le\Phi(d(a,a_0),d(a_0,b))
\leq\Phi(\alpha,\Phi(d(a_0,b_0),d(b_0,b)))
\leq\Phi(\alpha,\Phi(d(a_0,b_0),\beta)).
}
Therefore, with $\eps>\Phi(\alpha,\Phi(d(a_0,b_0),\beta))$, we obtain 
\Eq{*}{
 a\in \B(b,\eps)\subset\bigcup_{b'\in B}\B(b',\eps),\qquad
 b\in \B(a,\eps)\subset\bigcup_{a'\in A}\B(a',\eps).
}
and hence $D(A,B)\le\eps<+\infty$. It is clear that $D(A,B)=0$ if $A=B$.

Assume that $A,B\in\F(X)$ and $D(A,B)=0$. Let $a\in A$ fixed arbitrarily. Then, for every $n\in\N$, there is an element $b_n\in B$ such that $d(a,b_n)<1/n$. This means that the sequence $(b_n)$ tends to $a\in A$. Since $B$ is a closed set, therefore $a\in B$. However, $a\in A$ is arbitrary, so
$A\subset B$. The inclusion in the other direction can be proved similarly. That is, $A=B$, as it was desired.

The symmetry of $D$ is an immediate consequence of its definition.

Finally, for $A,B,C\in\F(X)$, we are going to show that
\Eq{ABC}{
 D(A,B) \leq \Phi\big(D(A,C),D(C,B)\big).
}
If $D(A,B)=0$, then there is nothing to prove. In the other case, let $\eps<D(A,B)$ be arbitrary. Then either $A\not\subseteq\bigcup_{b\in B}\B(b,\eps)$ or $B\not\subseteq\bigcup_{a\in A}\B(a,\eps)$. In the first case, there exists $a\in A$ such that $d(a,b)\geq\eps$ for all $b\in B$. Then, for all $c\in C$, we have 
\Eq{cb}{
\eps\leq d(a,b)\leq\Phi(d(a,c),d(c,b)).
}
We show that
\Eq{*}{
  \inf_{c\in C}d(a,c)\leq D(A,C).
}
If this were not true, then 
\Eq{*}{
  r:=\inf_{c\in C}d(a,c)> D(A,C).
}
Hence
\Eq{*}{
  a\in A\subseteq\bigcup_{c\in C}\B(c,r).
}
This inclusion says that $d(a,c)<r$ for some $c\in C$, which contradicts the definition of $r$. 

Thus, for every $n\in\N$ there exists $c_n\in C$ such that 
\Eq{*}{
  d(a,c_n)\leq D(A,C)+\frac1n.
}
Similarly, 
\Eq{*}{
  \inf_{b\in B}d(c_n,b)\leq D(C,B).
}
Therefore, there exists $b_n\in B$ such that 
\Eq{*}{
  d(c_n,b_n)\leq D(C,B)+\frac1n.
}
Thus, \eq{cb} implies
\Eq{*}{
\eps\leq \Phi(d(a,c_n),d(c_n,b_n))
 \leq \Phi\big(D(A,C)+\tfrac1n,D(C,B)+\tfrac1n\big).
}
Using that $\Phi$ upper semicontinuous, taking the limit $n\to\infty$ implies
\Eq{*}{
\eps \leq \Phi\big(D(A,C),D(C,B)\big).
}
Since $\eps<D(A,B)$ was arbitrary, we obtain that \eq{ABC} is valid.
This proves that $\Phi$ is a triangle function for $D$.
\end{proof}


One of the most important property of the fractal space is that it inherits the completeness of the base space. This is formulated by the following analogue of Blaschke's theorem which we describe in two versions.

\Thm{Blaschke1}{Let $(X,d)$ be a complete semimetric space with an upper semicontinuous triangle function $\Phi\colon\R_+^2\to\R_+$ for $d$. Then $(\F(X),D)$ is also complete.}

\begin{proof}
Let $(A_n)$ be an arbitrary Cauchy sequence in $(\F(X),D)$. We prove that $D(A_n,A)\to0$, where
\Eq{*}{
 A=\bigl\{x\in X\mid \exists(x_k):\,x_k\to x,\,x_k\in A_k \,\,(k\in\N)\bigr\}.
}

Let $\eps>0$ be arbitrary and construct a decreasing sequence $(\eps_k)$ of positive numbers converging to $0$ and satisfying the following properties:
\Eq{*}{
  \eps_1:=\eps,\qquad 
  \Phi(\eps_{k+1},\eps_{k+1})<\eps_k \qquad(k\in\N).
}
Due to the Cauchy property there exists $n_0\in\N$ such that $D(A_n,A_m)<\eps_2$ whenever $n,m>n_0$. We show that for every $n>n_0$ the following two inclusions hold:
\Eq{*}{
 \text{(i)} \qquad A\subset\bigcup_{y\in A_n}\B(y,\eps);\qquad
 \text{(ii)} \qquad A_n\subset\bigcup_{x\in A}\B(x,\eps).
}
First, let $x\in A$ be arbitrary. In this case, there is a sequence $(x_k)$ such that $x_k\in A_k$ and $x_k\to x$. Choose $k\in\N$ to satisfy $k>n_0$ and $d(x_k,x)<\eps_2$ simultaneously. Due to $D(A_k,A_n)<\eps_2$ there exists $y\in A_n$ such that $d(x_k,y)<\eps_2$. Therefore,
\Eq{*}{
 d(x,y)\le \Phi(d(x,x_k),d(x_k,y))
 \leq\Phi(\eps_2,\eps_2)<\eps_1=\eps.
}
Hence, $x\in \B(y,\eps)$, which proves the inclusion (i). 

To prove (ii), let $y\in A_n$ be arbitrary. Then, in view of the Cauchy property of the sequence $(A_k)$, there exists a strictly increasing sequence $(\ell_j)$ such that $\ell_1=n$, and for all $k>\ell_j$, the inequality
\Eq{kj}{
 D(A_{\ell_j},A_k)<\eps_{j+2} 
}
be valid. Now we construct a sequence $(x_k)$ in the following way. For $k<n$, let $x_k\in A_k$ be arbitrary. If $k=n$, then let $x_k=y$. Finally, for $k>\ell_1=n$, there exist $j\in\N$ such that $k\in\{\ell_j+1,\ldots,\ell_{j+1}\}$. Then, by \eq{kj}, there exists $x_k\in A_k$ such that $d(x_{\ell_j},x_k)<\eps_{j+2}$ be valid. We show that the sequence $(x_k)$ so constructed is a Cauchy sequence.

First, for $1\leq i\leq j$, we prove that
\Eq{*}{
  d(x_{\ell_i},x_{\ell_j})<\eps_{i+1}.
}
We are going to show, by induction with respect to $k$, that
this inequality holds for all $1\leq i\leq j\leq i+k$. If $k=0$, then $i=j$ and the inequality is trivial. Assume that the statement holds for some $k\geq0$ and let $1\leq i\leq j\leq i+k+1$. If $j<i+k+1$ (i.e., $j\leq i+k$), then the assertion follows from the inductive hypothesis. Thus, we may assume that $j=i+k+1$. Then, by the construction of the sequence $(\ell_i)$, $d(x_{\ell_i},x_{\ell_{i+1}})<\eps_{i+2}$ and by the inductive assumption, we have that $d(x_{\ell_{i+1}},x_{\ell_{i+k+1}})<\eps_{i+2}$. Therefore, using the triangle inequality, we obtain
\Eq{*}{
  d(x_{\ell_i},x_{\ell_j})=d(x_{\ell_i},x_{\ell_{i+k+1}})\leq\Phi(d(x_{\ell_i},x_{\ell_{i+1}}),d(x_{\ell_{i+1}},x_{\ell_{i+k+1}}))
  \leq\Phi(\eps_{i+2},\eps_{i+2})<\eps_{i+1}.
}
Let $l,m\in\N$ with $l,m>\ell_2$ be arbitrary. Then there exist $i,j\geq 2$ such that
\Eq{*}{
 l\in\{\ell_i+1,\ldots,\ell_{i+1}\},\qquad
 m\in\{\ell_j+1,\ldots,\ell_{j+1}\}.
}
We may assume that $l\le m$, then we also have that $i\le j$. Therefore,
\Eq{*}{
 d(x_l,x_m) 
 &\leq \Phi(d(x_l,x_{\ell_i}),d(x_{\ell_i},x_m))
 \leq \Phi(d(x_l,x_{\ell_i}),\Phi(d(x_{\ell_i},x_{\ell_j}),d(x_{\ell_j},x_m)))\\
 &\leq \Phi(\eps_{i+2},\Phi(\eps_{i+1},\eps_{j+2}))
 \leq \Phi(\eps_{i+2},\Phi(\eps_{i+1},\eps_{i+1}))
 \leq \Phi(\eps_{i+2},\eps_{i})
 \leq \Phi(\eps_{i},\eps_{i})<\eps_{i-1}.
}
Let $\eta>0$ be arbitrary. The sequence $(\eps_i)$ being a nullsequence,  there exists $i\geq2$ such that $\eps_{i-1}<\eta$ and then, for $m\geq l>\ell_i$, we have that $d(x_l,x_m)<\eps_{i-1}<\eta$. This proves that $(x_k)$ is a Cauchy sequence.

Due to the completeness of $X$, there exists $x\in X$ such that $x_k\to x$. Obviously, we have that $x\in A$. In view of the construction of the sequence $(x_{\ell_j})$, we have $d(y,x_{\ell_j})=d(x_{\ell_1},x_{\ell_j})<\eps_2$. Choose $j\geq1$ so that $d(x_{\ell_j},x)<\eps_2$. Thus,
\Eq{*}{
 d(y,x)=d(x_{\ell_1},x) \leq\Phi(d(x_{\ell_1},x_{\ell_j}),d(x_{\ell_j},x))
 \leq\Phi(\eps_2,\eps_2)<\eps_1=\eps.
 }
This shows that $y\in \B(x,\eps)$, whence we obtain that the inclusion $(ii)$ is also valid. The two inclusions imply that $D(A_n,A)\leq\eps$ whenever $n>n_0$. 

To complete the proof, it suffices to show that $A$ is a nonempty, bounded and closed set.

If $n>n_0$ then $(ii)$ is valid and $A_n$ is nonempty, thus $A$ must be nonempty, too. Similarly, $A$ is bounded because of $(i)$, since $A_n$ is bounded as well.
Indeed, if $x_1,x_2\in A$, then by $(i)$, there exist $y_1,y_2\in A_n$ such that
$d(x_i,y_i)<\eps$. Therefore,
\Eq{*}{
  d(x_1,x_2)\leq\Phi(d(x_1,y_1),d(y_1,x_2))
  \leq \Phi(d(x_1,y_1),\Phi(d(y_1,y_2),d(y_2,x_2)))
  \leq \Phi(\eps,\Phi(\diam(A_n),\eps)).
}
Hence, 
\Eq{*}{
  \diam(A)\leq \Phi(\eps,\Phi(\diam(A_n),\eps)),
}
which proves that $A$ is bounded.

To prove the closedness of $A$, let $z$ be an arbitrary element of the closure of $A$. Then, there exists a sequence $z_n\in A$ such that
\Eq{*}{
  d(z_n,z)<D(A_n,A)+\frac{1}{n}=:\delta_n.
}
Then $(\delta_n)$ is a null sequence which satisfies $D(A_n,A)<\delta_n$ for all $n\in\N$. Thus,
\Eq{*}{
  z_n\in A\subseteq\bigcup_{x\in A_n}\B(x,\delta_n)
}
and hence, there exists $x_n\in A_n$ such that $d(x_n,z_n)<\delta_n$. Then,
\Eq{*}{
  d(x_n,z)\leq\Phi(d(x_n,z_n),d(z_n,z))\leq\Phi(\delta_n,\delta_n).
}
By the upper semicontinuity of $\Phi$ at $(0,0)$ and $\Phi(0,0)=0$, it follows that the right hand side of the above inequality is a null sequence. Hence, $x_n\to z$, which shows that $z\in A$ and proves that $A$ is closed.
\end{proof}

\Thm{Blaschke2}{Let $(X,d)$ be a complete semimetric space with an upper semicontinuous triangle function $\Phi\colon\R_+^2\to\R_+$ for $d$. Then $(\K(X),D)$ is also complete.}

\begin{proof} In view of \thm{Blaschke1}, it suffices to show that if $(A_n)$ is a sequence of nonempty compact sets which converges to $A$ by the metric $D$, then $A$ is also compact.
In view of \thm{Haus}, it suffices to show that $A$ is totally bounded.

Let $\eps>0$ be arbitrary. Then there exists $\delta>0$ such that $\Phi(\delta,\delta)<\eps$. By the convergence $D(A_n,A)\to0$, there exists $n\in\N$ such that
\Eq{*}{
  A\subseteq\bigcup_{y\in A_n}\B(y,\delta).
}
On the other hand, by the total boundedness of $A_n$, there exist $p_1,\dots,p_m\in X$ such that 
\Eq{*}{
A_n\subseteq\bigcup_{i=1}^m\B(p_i,\delta).
}
We are going show that
\Eq{AA}{
A\subseteq\bigcup_{i=1}^m\B(p_i,\eps).
}
Indeed, if $a\in A$, then $a\in\B(y,\delta)$ for some $y\in A_n$. That is, $d(a,y)<\delta$. Furthermore, $y\in A_n$ implies that $y\in\B(p_i,\delta)$ for some $i\in\{1,\dots,n\}$, which means that $d(y,p_i)<\delta$. Therefore,
\Eq{*}{
  d(a,p_i)\leq\Phi(d(a,y),d(y,p_i))
  \leq\Phi(\delta,\delta)<\eps,
}
and hence $a\in \B(p_i,\eps)$ for some $i\in\{1,\dots,n\}$.
This proves \eq{AA} and shows that $A$ is totally bounded.
\end{proof}

\section{Fixed point theorems}

A mapping $\varphi:\R_+\to\R_+$ is called a \emph{comparison function} if it is increasing and $\lim_{n\to\infty}\varphi^n(t)=0$ for all $t>0$. It easily follows that $\varphi(t)<t$ holds for all $t>0$. Given a comparison function $\varphi$, a selfmap $T$ of a semimetric space $(X,d)$ is called a \emph{$\varphi$-contraction} if, for all $x,y\in X$,
\Eq{*}{
  d(T(x),T(y))\leq\varphi(d(x,y)).
}

The $\varphi$-contraction property of $T$ implies that, for all $x,y\in X$,
\Eq{*}{
  d(T(x),T(y))\leq\varphi(d(x,y))\leq d(x,y),
}
Thus, $T$ is Lipschitzian with a modulus $L=1$. According to \lem{cont}, this implies that $T$ is continuous.

In this setting, the following fixed point theorem was established by Bessenyei and Páles in \cite{BesPal17}.

\Thm{BesPal}{Let $(X,d)$ be a complete regular semimetric space and $\varphi$ be a comparison function and let $T\colon X\to X$ be $\varphi$-contraction. Then $T$ has  a unique fixed point $x_0\in X$ and, for all $x\in X$, the sequence $(x_k)$ defined by $x_1:=x$, $x_{k+1}:=T(x_k)$ converges to $x_0$.}

The following lemma allows us to establish error estimates concerning the iteration defined in the above theorem.

\Lem{est}{Let $(X,d)$ be a complete regular semimetric space with a triangle function $\Phi$ for $d$, let $\varphi$ be a comparison function, and let $T\colon X\to X$ be $\varphi$-contraction. Let $x_0$ be the unique fixed point of $T$. Then, for all $x\in X$,
\Eq{HK}{
  d(x,x_0) \le \psi(d(x,T(x))),
}
where $\psi:\R_+\to\overline{\R}_+$ is defined by
\Eq{psi}{
  \psi(t):=\sup\{s\geq0\mid s\leq \Phi(t,\varphi(s))\} \qquad (t\in\R_+).
}}

\begin{proof} Observe that $\psi$ is a nondecreasing function.

Let $x\in X$ be arbitrary. By the triangle inequality and by the $\varphi$-contractivity of $T$, we get
\Eq{*}{
  d(x,x_0)
  &\leq\Phi(d(x,T(x)),d(T(x),x_0))\\
  &=\Phi(d(x,T(x)),d(T(x),T(x_0)))
  \leq\Phi(d(x,T(x)),\varphi(d(x,x_0))).
}
Therefore,
\Eq{*}{
  d(x,x_0)\in \{s\geq0\mid s\leq \Phi(d(x,T(x)),\varphi(s))\},
}
which implies that
\Eq{*}{
  d(x,x_0)\leq\sup\{s\geq0\mid s\leq \Phi(d(x,T(x)),\varphi(s))\}
  =\psi(d(x,T(x))).
}
Thus the proof is complete.
\end{proof}

\Thm{error}{Let $(X,d)$ be a complete regular semimetric space with a triangle function $\Phi$ for $d$, let $\varphi$ be a comparison function and let $T\colon X\to X$  be a $\varphi$-contraction. Let $x_0\in X$ be the unique fixed point of $T$, let $x\in X$ and define the sequence $(x_k)$ recursively: $x_1:=x$, $x_{k+1}:=T(x_k)$. Then, for all $k\in\N$, the following inequalities hold:
\Eq{*}{
 d(x_k,x_0) &\le \psi\circ\varphi^{k-\ell}(d(x_k,x_{k+1})),\qquad \ell\in\{1,\dots,k\},\\
 d(x_k,x_0) &\le \varphi(d(x_{k-1},x_0)),
}
where $\psi:\R_+\to\overline{\R}_+$ is defined by \eq{psi}.}

If $\ell=1$ and $\ell=k$, then the first inequality is called the \emph{a priori} and \emph{a posteriori} error estimate, respectively. On the other hand, the second inequality is called the \emph{convergence speed} estimate.

\begin{proof} Let $\ell\in\{1,\dots,k\}$ and apply the inequality \eq{HK} of \lem{est} for $x:=x_k$. Then, using $(k-\ell)$ times the $\varphi$-contractivity of $T$ and the monotonicity of $\psi$, we get
\Eq{*}{
 d(x_k,x_0)&\le\psi(d(x_k,T(x_k)))
 =\psi(d(x_k,x_{k+1}))\\
 &=\psi(d(T^{k-\ell}(x_\ell),T^{k-\ell}(x_{\ell+1})))
 \leq \psi(\varphi^{k-\ell}(d(x_\ell,x_{\ell+1}))).
}
The second inequality is an immediate consequence of the $\varphi$-contraction property. Indeed, 
\Eq{*}{
   d(x_k, x_0)=d(T(x_{k-1}), T(x_0))\leq\varphi(d(x_{k-1}, x_0)),
}
which completes the proof.
\end{proof}

The following result from the paper \cite[Theorem 2]{BesPal17} establishes the stability of fixed points for a sequence of $\varphi$-contractions.

\Thm{BesPal2}{Let $(X,d)$ be a complete regular semimetric space, let $\varphi$ be a comparison function and let $T_k\colon X\to X$  be a sequence of $\varphi$-contractions which converges pointwise to a $\varphi$-contraction $T_0\colon X\to X$. For $k\in\N\cup\{0\}$, let $x_k$ denote the unique fixed point of $T_k$. Then $(x_k)$ converges to $x_0$.}

Based on \thm{BesPal} we can now present the first main result of this paper.

\Thm{Hutchinson}{Let $(X,d)$ be a complete semimetric space with an upper semicontinuous triangle function $\Phi\colon\R_+^2\to\R_+$ for $d$, let $\varphi:\R_+\to\R_+$ be an upper semicontinuous comparison function, and let $T_1,\ldots, T_n\colon X\to X$ be $\varphi$-contractions. Then there exists a unique fractal in $\K(X)$ with respect to the system $(T_1,\ldots,T_n)$.}

\begin{proof} Let the mapping $T:\K(X)\to\K(X)$ be defined in the following way:
\Eq{T}{
 T(H)=\bigcup_{k=1}^nT_k(H).
}
If $H$ is a compact set, then, by the continuity of each $T_k$, the set $T_k(H)$ is also compact and hence $T(H)$ is compact, which proves that $T$ maps $\K(X)$ into itself.

We are going to show that $T$ is a $\varphi$-contraction on $(\K(X),D)$. Let $A,B\in\K(X)$ and $\varphi(D(A,B))<\eps$. By the upper semicontinuity of $\varphi$, there exists $D(A,B)<\delta$ such that $\varphi(\delta)<\eps$. If $a\in A$ is arbitrary, then there exists $b\in B$ such that $d(a,b)<\delta$. Therefore, for every $k\in\{1,\ldots,n\}$, we have
\Eq{*}{
 d(T_k(a),T_k(b))\le \varphi(d(a,b))\leq\varphi(\delta)<\eps
}
and hence,
\Eq{*}{
 T_k(a)\in \B(T_k(b),\eps)\subset\bigcup_{y\in T(B)}\B(y,\eps).
}
Since $a\in A$ and $k\in\{1,\ldots,n\}$ were arbitrary, we obtain that
\Eq{*}{
 T(A)\subset\bigcup_{y\in T(B)}\B(y,\eps).
}
It can be proved similarly that
\Eq{*}{
 T(B)\subset\bigcup_{x\in T(A)}\B(x,\eps).
}
Therefore, $D\bigl(T(A),T(B)\bigr)\le\eps$. Upon taking the right-hand-side limit $\eps\to \varphi(D(A,B))$, it follows that $D\bigl(T(A),T(B)\bigr)\leq \varphi(D(A,B))$, which shows that $T$ is a $\varphi$-contraction on $(\K(X),D)$.

In view of \thm{Blaschke2}, we have that $(\K(X),D)$ is a complete semimetric space. Thus, applying the fixed point theorem \thm{BesPal}, we obtain that $T$ has a unique fixed point $H_0\in\K(X)$, which is the fractal with respect to the function system $(T_1,\ldots,T_n)$.
\end{proof}

\Thm{festimates}{Let $(X,d)$ be a complete regular semimetric space with an upper semicontinuous triangle function $\Phi\colon\R_+^2\to\R_+$ for $d$, let $\varphi:\R_+\to\R_+$ be an upper semicontinuous comparison function, and let $T_1,\ldots T_n\colon X\to X$ be $\varphi$-contractions. Let $H\in\K(X)$ and let $H_0$ be the unique fractal in $\K(X)$ with respect to the system $(T_1,\ldots, T_n)$. Define the sequence $(H_k)$ recursively: $H_1:=H$, $H_{k+1}:=T(H_k)$. Then, for all $k\in\N$ the following inequalities hold:
\Eq{*}{
 D(H_k,H_0) &\le \psi\circ\varphi^{k-\ell}(D(H_k,H_{k+1})),\qquad \ell\in\{1,\dots,k\}\\
 D(H_k,H_0) &\le \varphi(D(H_{k-1},H_0)),
}
where $\psi:\R_+\to\overline{\R}_+$ is defined by \eq{psi}.}

\begin{proof} The statement is direct consequence of \thm{error} if we apply it to the semimetric space $(\K(X),D)$ and the $\varphi$-contraction $T:\K(X)\to\K(X)$ defined by \eq{T}.
\end{proof}

\Thm{stability}{Let $(X,d)$ be a complete regular semimetric space with an upper semicontinuous triangle function $\Phi\colon\R_+^2\to\R_+$ for $d$ and let $\varphi:\R_+\to\R_+$ be an upper semicontinuous comparison function. Let $(T_{k,1}),\dots,(T_{k,n})$ be sequences of $\varphi$-contractions that converge pointwise to some $\varphi$-contractions $T_{0,1},\dots,T_{0,n}$, respectively. For $k\in\N\cup\{0\}$, define the mapping $\mathbb{T}_k\colon\K(X)\to\K(X)$ by
\Eq{*}{
 \mathbb{T}_k(H)=\bigcup_{j=1}^nT_{k,j}(H)
}
and denote the unique fractal in $\K(X)$ with respect to the system $(T_{k,1},\ldots, T_{k,n})$ by $H_k$. Then the sequence $(H_k)$ converges to $H_0$ with respect to the semimetric $D$.
}

\begin{proof} By the pointwise convergence, for all $j\in\{1,\dots,n\}$, we have that $D(T_{k,j}(H), T_{0,j}(H))\longrightarrow 0$. Applying \lem{H_k}, it follows that $D(\mathbb{T}_k(H),\mathbb{T}_0(H))\longrightarrow0$, that is, $\mathbb{T}_k$ converges to $\mathbb{T}_0$ pointwise. As we have proved it in the proof of \thm{Hutchinson}, $\mathbb{T}_k$ is a $\varphi$-contraction on $\K(X),D)$. Therefore, the statement is a direct consequence of \thm{BesPal2}. 
\end{proof}


\end{document}